\documentclass[10pt,a4paper]{article}
\usepackage{graphicx}
\usepackage{amsmath}
\usepackage{amsfonts}
\usepackage{amssymb}
\usepackage{amsthm}
\usepackage{color}
\begin{document}

\def\R{\mathbb{R}}
\def\N{\mathbb{N}}
\def\H{\mathcal{H}}
\def\d{\textrm{div}}
\def\v{\textbf{v}}
\def\I{\hat{I}}
\def\B{\hat{B}}
\def\x{\hat{x}}
\def\y{\hat{y}}
\def\p{\hat{\phi}}
\def\r{\hat{r}}
\def\w{\textbf{w}}
\def\u{\textbf{u}}
\def\K{\mathcal{K}}
\def\Reg{\textrm{Reg}}
\def\s{\textrm{Sing}}
\def\sgn{\textrm{sgn}}

\newtheorem{defin}{Definition}[section]
\newtheorem{lem}{Lemma}[section]
\newtheorem{rem}{Remark}[section]
\newtheorem{cor}{Corollary}[section]
\newtheorem{thm}{Theorem}[section]
\newtheorem{prop}{Proposition}[section]
\newtheorem{definition}{Definition}[section]
\newtheorem{con}{Conjecture}[section]
\newtheorem{Main}{Main Result}

\title{Short note on energy maximization property of the first eigenfunction of the Laplacian}

\author{Hayk Mikayelyan \thanks{Mathematical Sciences, University of Nottingham Ningbo, 199 Taikang East Road, Ningbo 315100, PR China \hfill Hayk.Mikayelyan@nottingham.edu.cn  } }

\maketitle

\begin{abstract}
We consider the Dirichlet-energy maximization problem of the solution 
$u_f$ of (\ref{main}), among all functions $f\in  L^2(D)$, such that 
$\|f\|_2= 1$. We show that the two maximizers are the first eigenfunctions of the Laplacian with Dirichlet boundary condition $f=\pm u_1$.
\end{abstract}

%\tableofcontents

%%%%%%%%%%%%%%%%%%%%%%%%%%%%%%%%%%%%%%%%%%%%%%%%%%%%%%%%%%%%%%%%%%%%%%%%%%%%%%%%%%%%%%%%%%%%%%%%%%%%%%%%%%%%%%%%%%%%%%%%%%%%%%%
%%%%%%%%%%%%%%%%%%%%%%%%%%%%%%%%%%%%%

One of the classical problems in rearrangements theory is the minimization/maximization of the
functional
\begin{equation}\label{Dir-functional}
\Phi(f)=\int_D |\nabla u_f|^2 dx,
\end{equation}
where $u_f$ is the unique solution of the Dirichlet boundary value problem
\begin{equation}\label{main}
\begin{cases} -\Delta u_f = f  & \mbox{in } D,
 \\
 u_f=0  & \mbox{on } \partial D,  \end{cases}
\end{equation}
and $f$ belongs to the rearrangement class.

Let us change the problem and maximize the functional (\ref{Dir-functional}) over the unit sphere in $L^2$, i.e.,
$$
f\in \mathcal{S}=\{f\in L^2(D)\,:\, \|f\|_2=1  \}.
$$

%%%%%%%%%%%%%%%%%%%%%%%%%%%%%%%%%%%%%%%%%%%%%%%%%%%%%%%%%%%%%%%%%%%%%%%%%%%%%%%%%%%%%%%%%%%%%%%%%%%%%%%%%%%%%%%%%%%%%%%%%%%%%%%
%%%%%%%%%%%%%%%%%%%%%%%%%%%%%%%%%%%%%%%%%%%%%%%%%%%%%%%%%%%%%%%%%%%%%%%%%%%%%%%%%%%%%%%%%%%%%%%%%%%%%%%%%%%%%%%%%%%%%%%%%%%%%%%

First we relax the problem and consider the maximization of $\Phi(f)$ over the unit ball
$$
\mathcal{B}=\{f\in L^2(D)\,:\, \|f\|_2\leq 1  \}
$$
in $L^2(D)$. The existence of the maximizer(s) $\hat{f}\in \mathcal{S}$ follows from weak closedness and convexity of $\mathcal{B}$, and strict convexity and weak continuity of $\Phi$.

Passing to the limit in the extremality condition
$$
t^{-1}\left[\Phi(\hat{f}+t(f-\hat{f}))-\Phi(\hat{f})\right]\leq 0
$$
we obtain
$$
\langle \Phi'(\hat{f}),f-\hat{f} \rangle \leq 0.
$$
Observe that $\Phi'(\hat{f})$ can be associated with $2\hat{u}:=2u_{\hat{f}}$,
\begin{multline}
\epsilon^{-1}\int_D |\nabla u_{f+\epsilon h}|^2-|\nabla u|^2 dx=\\  \int_D\nabla(u_{f+\epsilon h}+u_f)\cdot
\nabla u_h dx\to_{\epsilon\to 0 } \int_D 2u_fhdx.
\end{multline}

Thus,
$$
\int_D f\hat{u} \leq \hat{f}\hat{u}dx
$$
for all $f\in\mathcal{B}$.

Applying Cauchy-Schwartz inequality we see that
$$
\int_D \hat{f}\hat{u}dx\leq \|\hat{u}\|_2\|\hat{f}\|_2\leq\|\hat{u}\|_2
$$
and the equality holds if and only if $\hat{f}=\lambda\hat{u}$.
Thus, $\hat{f}$ must coincide with one of the eigenfunctions $u_k\in \mathcal{S}$ of the eigenvalue problem with Dirichlet boundary conditions (see \cite{Henrot:book})
\begin{equation}\label{eigen}
\begin{cases} -\Delta u = \lambda u   & \mbox{in } D,
 \\
 u=0  & \mbox{on } \partial D ,  \end{cases}
\end{equation}
and $\hat{u}=\lambda^{-1}_k(D)u_k$, where $\lambda_k$ is the $k$th eigenvalue.

The maximization of $\Phi$ over $\mathcal{B}$ reduces now to the maximization of
$$
\Phi(u_k)=\lambda^{-2}_k(D) \int_D |\nabla u_k|^2=\frac{1}{\lambda_k(D)}
$$
over $k$, which happens when $k=1$.

We have proven the following theorem

\begin{thm}
The solutions of the maximization problem
$$
\max_{\|f\|_2\leq 1} \int_D |\nabla u_f|^2dx,
$$
where the function $u_f$ is the solution of the equation (\ref{main}), are the
first eigenfunctions $\hat{f}=\pm u_1\in \mathcal{S}$ of the Dirichlet eigenvalue problem (\ref{eigen}),
$$
\Phi(\hat{f})=\frac{1}{\lambda_1(D)},
$$
and $\hat{u}=\pm \lambda^{-1}_1(D)u_1$.
\end{thm}
%%%%%%%%%%%%%%%%%%%%%%%%%%%%%%%%%%%%%%%%%%%%%%%%%%%%%%%%%%%%%%%%%%%%%%%%%%%%%%%%%%%%%%%%%%%%%%%%%%%%%%%%%%%%%%%%%%%%%%%%%%%%%%%
%%%%%%%%%%%%%%%%%%%%%%%%%%%%%%%%%%%%%%%%%%%%%%%%%%%%%%%%%%%%%%%%%%%%%%%%%%%%%%%%%%%%%%%%%%%%%%%%%%%%%%%%%%%%%%%%%%%%%%%%%%%%%%%

%%%%%%%%%%%%%%%%%%%%%%%%%%%%%%%%%%%%%%%%%%%%%%%%%%%%%%%%%%%%%%%%%%%%%%%%%%%%%%%%%%%%%%%%%%%%%%%%%%%%%%%%%%%%%%%%%%%%%%%%%%%%%%%
%%%%%%%%%%%%%%%%%%%%%%%%%%%%%%%%%%%%%%%%%%%%%%%%%%%%%%%%%%%%%%%%%%%%%%%%%%%%%%%%%%%%%%%%%%%%%%%%%%%%%%%%%%%%%%%%%%%%%%%%%%%%%%%
\begin{rem}
To see that the minimization of $\Phi$ over $\mathcal{S}$ does not have a solution we need only to take a
sequence $f_k\in \mathcal{S}$, such that $f_k \rightharpoonup 0$ and see that $\Phi(\hat{u}_k)\to 0$. Observe that $f=0$ is the minimizer of $\Phi$ over $\mathcal{B}$.
\end{rem}
%%%%%%%%%%%%%%%%%%%%%%%%%%%%%%%%%%%%%%%%%%%%%%%%%%%%%%%%%%%%%%%%%%%%%%%%%%%%%%%%%%%%%%%%%%%%%%%%%%%%%%%%%%%%%%%%%%%%%%%%%%%%%%%
%%%%%%%%%%%%%%%%%%%%%%%%%%%%%%%%%%%%%%%%%%%%%%%%%%%%%%%%%%%%%%%%%%%%%%%%%%%%%%%%%%%%%%%%%%%%%%%%%%%%%%%%%%%%%%%%%%%%%%%%%%%%%%%

%%%%%%%%%%%%%%%%%%%%%%%%%%%%%%%%%%%%%%%%%%%%%%%%%%%%%%%%%%%%%%%%%%%%%%%%%%%%%%%%%%%%%%%%%%%%%%%%%%%%%%%%%%%%%%%%%%%%%%%%%%%%%%%
%%%%%%%%%%%%%%%%%%%%%%%%%%%%%%%%%%%%%%%%%%%%%%%%%%%%%%%%%%%%%%%%%%%%%%%%%%%%%%%%%%%%%%%%%%%%%%%%%%%%%%%%%%%%%%%%%%%%%%%%%%%%%%%

\bibliographystyle{plain}
\bibliography{rearr}

\begin{thebibliography}{1}

\bibitem{Henrot:book}
Antoine Henrot.
\newblock {\em Extremum problems for eigenvalues of elliptic operators}.
\newblock Frontiers in Mathematics. Birkh\"auser Verlag, Basel, 2006.

\end{thebibliography}

\end{document}